\documentclass[12pt,oneside]{amsart}

      \usepackage{amssymb}
      \usepackage{rotating}

\usepackage{amsfonts}
\usepackage{dsfont}
\usepackage{amsxtra}
\usepackage{graphicx}
\DeclareGraphicsExtensions{.pdf,.png,.jpg}

\usepackage{index}
 \newtheorem{thm}{Theorem}[section]
 \newtheorem{cor}[thm]{Corollary}
 \newtheorem{lemma}[thm]{Lemma}
 \newtheorem{prop}[thm]{Proposition}
 \theoremstyle{definition}
 \newtheorem{defn}[thm]{Definition}
 \theoremstyle{remark}
 \newtheorem{remark}[thm]{Remark}
 \numberwithin{equation}{section}
 \theoremstyle{note}
 
 \newtheorem{exam}[thm]{Example}
 
      \makeatletter
      \def\@setcopyright{}
      \def\serieslogo@{}
      \makeatother

 \catcode`\á=\active \def á{\'a}
 \catcode`\Á=\active \def Á{\'A}
  \catcode`\ó=\active \def ó{\'o}
  \catcode`\é=\active \def é{\'e}
  \catcode`\ú=\active \def ú{\'u}
  \catcode`\í=\active \def í{\'{\i}}
  \catcode`\ñ=\active \def ñ{\~n}
  \catcode`\Ñ=\active \def Ñ{\~N}
  \catcode`\¿=\active \def ¿{?`}
  \catcode`\º=\active \def º{$^{\underline{o}}$}
  \catcode`\ª=\active \def ª{$^{\underline{a}}$}
  \catcode`\¡=\active \def ¡{!`}
  \catcode`\â=\active \def â{\^{a}}
  \catcode`\ê=\active \def ê{\^{e}}
  \catcode`\î=\active \def î{\^{\i}}
  \catcode`\ô=\active \def ô{\^{o}}
  \catcode`\û=\active \def û{\^{u}}
  \catcode`\ç=\active \def ç{\c{c}}
  \catcode`\ü=\active \def ü{\"{u}}
  \catcode`\ö=\active \def ö{\"{o}}
\newenvironment{dem}{\noindent{\it Proof}.}{\qed}

\newcommand{\ff}{\mathbb{F}}

\newcommand{\oo}{\mathcal{O}}
\newcommand{\mm}{\mathfrak{m}}
\newcommand{\ago}{\mathfrak{a}}

\newcommand{\ooo}{\overline{\mathcal{O}}}

\begin{document}

% First we specify the top matter (author, title, etc).
%
% Note: All of the top matter items are optional and can be omitted.
% But you will probably want to specify at least the author and title
% and perhaps an abstract.

   % author information

   % first author 
   
%  \author{F\'elix Delgado de la Mata}
%  \address{Departamento de Algebra, Geometria y Topologia, Universidad de Valladolid. Prado de la Magdalena, s/n E-47005 Valladolid, Spain}
%  \email{fdelgado@agt.uva.es}
   
%   \author{Evgeny Gorsky}

   % the address where the research was carried out
%   \address{USA}

   % current address, usually not needed because it is the same as the
   % regular address
   %\curraddr{Department of Mathematics, Pennsylvania State University,
   %  University Park, State College PA 16802}

 %  \email{usa@usa.usa}

\author{Julio Jos\'e Moyano-Fern\'andez}
\address{Institut f\"ur Mathematik, Universit\"at Osnabr\"uck. Albrechtstra\ss e 28a, D-49076 Osnabr\"uck, Germany}
\email{jmoyanof@uni-osnabrueck.de}

   % title

\title[On the coefficients of the St\"ohr Zeta Function]
 {On the coefficients of the \\ St\"ohr Zeta Function}

   % Note that the short title for running heads goes in square
   % brackets.  This is optional.  The long title goes in curly
   % braces.  In the long title, line breaks are indicated by \\.

   % abstract (optional)
   \begin{abstract}
Let $\oo$ be a one-dimensional Cohen-Macaulay local ring having a finite field as a coefficient field.
 The aim of this work is to extend the explicit
computations of the St\"ohr Zeta Function of $\oo$
for one and two branches to an arbitrary number of them, obtaining
in this general case an upper bound for the coefficients
of the zeta function, instead of an equality. The calculations are
based on the use of the value semigroup of a curve singularity and a
suitable classification of the maximal points of the semigroup.
   \end{abstract}

   % AMS subject classifications (used in AMS journals)
   \subjclass{Primary 14H20; Secondary 32S10, 11G20}

   % AMS keywords (used in AMS journals)
   \keywords{Curve singularity, zeta function, Poincar\'e
series, value semigroup}

   % acknowledge support, etc
   \thanks{The author was partially supported by the Spanish Government grant ``Ministerio de Educaci\'on y Ciencia (MEC) 
MTM2007-64704", in cooperation with the European Union in the framework of the founds ``FEDER'',  and by the Deutsche Forschungsgemeinschaft (DFG)}
  % \thanks{We would like to thank our colleagues for their helpful
   %  criticism.?}

\maketitle

\section{Introduction}

Zeta functions on curve singularities over finite fields were introduced by Galkin \cite{galkin} and Green \cite{green} as attempts to reproduce the theory already known for the smooth case. A generalisation of these due to St\"ohr focused the description on techniques rather algebraic (cf. \cite{stohr}, \cite{stohr2}). The new approach allowed Delgado and the author to connect the St\"ohr zeta function with the Poincar\'e series defined on certain filtrations associated with curve singularities, building a bridge between the two viewpoints (see \cite{demo}).
\medskip

The coefficients of the zeta function are determined by the structure of the value semigroup associated to the singularity so that a good understanding of them contributes to a more precise knowledge of the singularity.
\medskip

The purpose of this paper is precisely to use a combinatorial description of the value semigroup for a better understanding of the coefficients of the St\"ohr zeta function. The first part of the work is devoted to present and summarise the state of the art concerning St\"ohr zeta function. In the second part we introduce a new combinatorial description of the value semigroup in order to obtain upper bounds for the coefficients of the St\"ohr zeta functions.
\medskip

\section{The value semigroup} \label{semigrupo:valores}

Let $\oo$ be a one-dimensional Cohen-Macaulay  local ring having a finite field
$\ff := \ff_q$ as a coefficient field. Let $\mm$ be its maximal
ideal. The extension degree $\rho$ of the residue field $\oo / \mm$ over
$\ff$ is finite (in fact, a power of $q$).  Denote by $K$ the
ring of fractions of $\oo$ and let $\ooo$ be the integral closure
of $\oo$ in $K$. Let $\ooo$ assume to be an $\oo$-module of
finite length.
\medskip

From now on, we distinguish between \emph{fractional} and
\emph{integral} (or proper) ideals, which are the usual ones. A
fractional ideal $\ago$ of $\oo$  is a submodule of $K$
such that there exists $z \in \oo$, $z \ne 0$ with $z \ago
\subseteq \oo$.
\medskip

We know that $K = \mathrm{Quot} (\oo) = \mathrm{Quot} (\ooo)$. So,
as it is shown in \cite{kiyek}, Prop. (2.10), there is a finite
number of \emph{Manis valuations} $v_1,
\ldots, v_r$ with associated valuation rings $V_1, \ldots, V_r$
such that $\ooo = V_1 \cap \ldots \cap V_r$. (Recall that a Manis
valuation is a surjective valuation defined over a ring in which
all prime ideal containing the Jacobson ideal is maximal and being
its own ring of fractions. This slightly extended concept of valuation is needed in
our context, since the ring $\oo$ is not a domain, and thus $K$ is
not a field. Further details in \cite{kiyek}, \S 2.2.) For each $i=1, \ldots
, r$, if $\mm (V_i)$ denotes the maximal ideal of $V_i$, then $\mm
(V_i) \cap \ooo = \mm_i$ is a regular maximal ideal (i.e.,
containing regular elements) of $\ooo$. These ideals are also
principal, whose generators are denoted by $\pi_i$. It is easy to see that the ideals of $\ooo$ are of the form $\mm^{\underline{n}} =
\mm_1^{n_1} \cdot \ldots \cdot \mm_r^{n_r}$, if
$\underline{n}=(n_1, \ldots , n_r) \in \mathbb{Z}^r$.
\medskip

Let $d_i := \dim_{\ff} \left ( \ooo / \mm_i \right )$, i.e., the
extension degree of the residue field of $\mm_i$ over $\ff$. We
have
\[
\dim_{\ff} \left ( \ooo / \mm^{\underline{n}}  \right ) = n_1 d_1
+ \ldots + n_r d_r := \underline{n} \cdot \underline{d},
\]
whenever $n_i \ge 0$ for all $1 \le i \le r$.
\medskip

The fractional ideals $\mathfrak{a}$ of $\oo$ are of the form
$\pi^{-\underline{n}} \mathfrak{b}$, with $\pi^{\underline{n}} :=
\pi_1^{n_1} \cdot \ldots \cdot \pi_r^{n_r}$ for $\underline{n} =
(n_1, \ldots , n_r) \in \mathbb{Z}^r$ and where $\mathfrak{b}$ is
an ideal of $\oo$ satisfying $\mathfrak{b} \cdot \overline{\oo} =
\overline{\oo}$. Indeed, if $\mathfrak{a}$ is a fractional ideal of $\oo$, then $z
\mathfrak{a}$ is a fractional ideal of $\oo$ for all $z \in
K^{\ast}:= K \setminus \{ 0 \}$. Moreover, $\pi^{\underline{m}}
\mathfrak{a}$ is a fractional ideal of $\oo$ for all
$\underline{m} \in \mathbb{Z}^r$. Let $-\underline{m} = (v_1(a),
\ldots, v_r (a))$ for some $a \in \mathfrak{a}$. If $\underline{m}
\ne \underline{0}$, then
\begin{itemize}
    \item[a)] $\pi^{\underline{m}} \mathfrak{a} = \mathfrak{b}$ is
    a fractional ideal of $\oo$.
    \item[b)] The set of units of $\ooo$ is 
    \[
    U_{\ooo}:= \{ z \in K \mid z \ooo = \ooo \}  = \{z \in K \mid v_i
(z)=0, ~ ~ 1 \le i \le r \}.
\]
    \item[c)] $\left ( v_1 \left ( \pi^{m_1} a \right ), \ldots, v_r \left ( \pi^{m_r} a \right ) \right ) =
    0$, hence $\pi^{-\underline{m}}a \in U_{\ooo}$.
\end{itemize}
Therefore $\mathfrak{a} = \pi^{-\underline{m}} \mathfrak{b}$.
\medskip

Finally, we do the assumption that the base field $\ff$ has a
cardinality equal to or greater than $r$.
\medskip

We present now an important invariant of (plane) curve
singularities, key to search common points between zeta functions
and singularity theory: the value semigroup
associated with a (plane) curve singularity, which is defined as
\[
S(\oo):= \left \{(v_1(z), \ldots, v_r(z)) \mid z \in
\oo \setminus \{0\} \right \}.
\]

We denote it
by $S$ whenever there is no risk of confusion. The semigroup $S$
is an additive sub-semigroup of $\mathbb{Z}_{+}^r$ which satisfies
the following properties (cf. \cite{de}):

\begin{itemize}
    \item[(S1)] $\underline{0} \in S$ and there exists $\min \left \{ S \setminus \{ \underline{0} \} \right
    \}\in S$.
    \item[(S2)] For $\underline{n},\underline{m} \in S$, the vector 
    \[
    \inf(\underline{n},\underline{m})=\left( \min \{n_1,m_1 \}, \ldots, \min \{n_r,m_r \} \right)
    \]
     belongs to
    $S$.
    \item[(S3)] $S$ has a conductor, i.e., there exists an element
    $\underline{\delta} \in S$ such that if $\underline{n} \in \mathbb{Z}_{+}^r$ with $\underline{n}
    \ge \underline{\delta}$, then $\underline{n} \in S$ and $\underline{\delta}$ is minimal with this property.
    \item[(S4)] (Pivotage) For $\underline{n},\underline{m} \in S$ such that there exists $i_0
    \in I=\{1, \ldots, r \}$ with $n_{i_0}=m_{i_0}$, then
    there is $\underline{\beta} \in S$ satisfying $\beta_k \ge \min \{ n_k, m_k \}$
    for any $k \in I$, $\beta_t = \min \{n_t,m_t \}$ if $n_t
    \ne m_t$ and $\beta_{i_0} > n_{i_0}=m_{i_0}$.
\end{itemize}

%Todas estas propiedades son de demostración sencilla. Para
%detalles complementarios, consultar \cite{de}.

\section{Maximal points of the value semigroup}

The classification of the elements in the semigroup of values will
be a useful tool to get explicit formulae of our zeta function.
\medskip

Let $\underline{n}=(n_1, \ldots, n_r) \in \mathbb{Z}_{+}^r$, $I=
\{ 1, \ldots, r \}$, $J \subset I$ with $J=\{ i_1, \ldots, i_h \}$
and $k \in I \setminus J$. We consider 
\begin{itemize}
    \item the set $\overline{\Delta}_J (\underline{n})$, which is defined to be
    \[
     \left \{ \underline{\beta} \in \mathbb{Z}_{+}^r \mid \beta_i = n_i ~
    \forall i \in J, ~ \beta_t > n_t ~ \forall t \notin J \right
    \}.
    \]
    In particular, for $i \in I$ we write
    \[
\overline{\Delta}_i (\underline{n})=\left \{ \underline{\beta} \in
\mathbb{Z}_{+}^r \mid \beta_i = n_i, ~ \beta_t > n_t ~ \forall t
\ne i \right \}.
    \]
By restricting to elements of the semigroup we employ the following notation:
    \[
\Delta_J(\underline{n})=\overline{\Delta}_J (\underline{n}) \cap
S.
    \]

    \item  the set $\overline{\Delta}_J^k(\underline{n})$, consisting of all $ \underline{\beta}
\in \mathbb{Z}_{+}^r$ such that
    \[
 \beta_i = n_i ~ \forall i \in J;~
\beta_d \ge n_d ~ \forall d \notin J, ~ d \le k; 
\beta_t > n_t ~
\forall t > k,~ t \notin J:
    \]
as before we write
    \[
\Delta_J^k (\underline{n})=\overline{\Delta}_J^{k} (\underline{n})
\cap S.
    \]
\end{itemize}
\medskip

%Procedemos a establecer una clasificación de los puntos de
%$\mathbb{Z}_{+}^r$ en relación al semigrupo de valores $S$ y su
%distribución en el espacio discreto $r$-dimensional.
%\medskip

The elements of the semigroup $S$ will be called points of the semigroup, the elements of $\mathbb{Z}_{+}^r \setminus S$ are called gaps of $S$. Among the points $\underline{n} \in S$, those satisfying $\Delta (\underline{n}) = \varnothing$ will be called maximal points of the semigroup.
\begin{defn}
A maximal point of $S$ such that $\Delta_J (\underline{n}) = \varnothing$ for all $J \subsetneq I$ is called absolute. A maximal point such that $\Delta_J (\underline{n}) \notin \varnothing$ for all $J \subsetneq I$ with $\sharp (J) \ge 2$ is called relative.
\end{defn}

The following definition generalises the concepts of relative and  absolute maximal points of $S$:

\begin{defn}
Let $x_2,x_3, \ldots, x_{r-1} \in \{0,1\}$. A maximal point $\underline{n} \in S$ is said to be of kind 
$x_2, \ldots , x_{r-1}$
 if
\[
x_{\sharp J}=\left\{%
\begin{array}{ll}
    0, & \hbox{if}~ \Delta_J (\underline{n})=\emptyset \mathrm{~for~all~} J \subsetneq  \mathrm{~with~} \sharp (J)=i; \\
    1, & \hbox{if}~ \Delta_J (\underline{n}) \ne \emptyset \mathrm{~for~all~} J \subsetneq  \mathrm{~with~} \sharp (J)=i; \\
\end{array}%
\right.
\]
\end{defn}

\begin{remark}
According to this definition, the maximal points of kind $0, \ldots,
0$ are the absolute maximal points, and those of kind $1,
\ldots, 1$ are relative ones. 
\end{remark}

\begin{exam}
In the case of two branches, the concepts of absolute maximal and
relative maximal coincide. For three branches, the set of maximal
points distinguishes only these two kinds. The case of four
branches needs a more precise study: a point $\underline{n} \in S
\subset \mathbb{Z}_{+}^4$ is maximal if $\Delta_i (\underline{n})
= \emptyset$ for every $i \in I=\{1,2,3,4 \}$. It will be:
\begin{itemize}
    \item[-] absolute maximal if moreover $\Delta_{ij}(\underline{n}) =
    \emptyset = \Delta_{ijk}(\underline{n})$.
    \item[-] relative maximal if $\Delta_{ij}(\underline{n}) \ne \emptyset$
    and $\Delta_{ijk}(\underline{n}) \ne \emptyset$.
\end{itemize}
Furthermore, these points are maximal, too:
\begin{itemize}
    \item[-] maximal points with $\Delta_{ij}(\underline{n})=\emptyset$
    but $\Delta_{ijk}(\underline{n}) \ne \emptyset$, i.e., maximal of
    kind $0,1$.
    \item[-] maximal points with $\Delta_{ij}(\underline{n}) \ne \emptyset$
    but $\Delta_{ijk} (\underline{n}) = \emptyset$, that is, those
    of kind $1,0$.
\end{itemize}
\end{exam}

Two basic properties are summarised in the next lemmas (see \cite[Lemma 1.3.4, Lemma 1.3.5, pp. 101-103]{de}):

\begin{lemma}
Let $\underline{n} \in \mathbb{Z}_{+}^r$ and assume that
$\underline{n} \notin S$. Then there exists one index $i \in I$ so
that $\Delta^r_i (\underline{n}) = \emptyset$.
\end{lemma}
\begin{lemma}
Let $\underline{n} \in \mathbb{Z}_{+}^r$ satisfying
\begin{enumerate}
    \item There is $i \in I$ with $\Delta_i (\underline{n}) = \emptyset$.
    \item For every $j \in I$, $j \ne i$, we have
    $\Delta_{ij}(\underline{n})\ne \emptyset$.
\end{enumerate}
Then $\underline{n} \in S$ and $\underline{n}$ is a relative
maximal point of $S$.
\end{lemma}
\medskip

Let $J=\{i_1, \ldots, i_h \} \subset I=\{1, \ldots,r \}$. Denoting
by $\mathrm{pr}_J$ the projection over the indexes $J$ of
$\mathbb{Z}_{+}^{r}$ in $\mathbb{Z}_{+}^{\sharp J}$. We set
\[
S_{i_1 , \ldots, i_h}=S_J:=\mathrm{pr}_J (S)= \{\mathrm{pr}_J
(\alpha) \mid \alpha \in S \}.
\]
That is, the semigroup $S_J$ corresponds to the projection of
$S_I$ on the branches indexed by $J$.

\begin{remark}
Let $\underline{n}$ be a gap of $S$. It could be caused either by a gap in any of the projection semigroups, or by the existence of a maximal point 
$\underline{\gamma}$ with $\gamma_i \le n_i$ for all $i \in \{1,
\ldots,r\}$.
\end{remark}

\section{The St\"ohr Zeta Function}

%Set-up: wer ist $\oo$, wer ist $K$, wer ist $\ff$, wer sind die
%$v_i$. (Am wenigsten).
%\medskip

For (integral) ideals of $\oo$, Galkin introduced in \cite{galkin}
the zeta function $\zeta_{\oo} (s):= \sum_{\mathfrak{a} \subseteq
\oo} \sharp \left (\oo / \mathfrak{a} \right )^{-s}; ~ \mathrm{Re}
(s) >0$. By considering fractional ideals $\ago \supset \oo$ instead,
St\"ohr defines in \cite{stohr2} the series
\[
\zeta (\oo, s):= \sum_{\mathfrak{a} \supseteq \oo} \sharp \left
(\mathfrak{a}/ \oo \right )^{-s}; ~ \mathrm{Re} (s) >0,
\]
which is called the \textbf{St\"ohr zeta function} of $\oo$. To
compare both zeta functions, the following Dirichlet series is
introduced:
\[
\zeta (\mathfrak{d},s):= \sum_{\mathfrak{a} \supseteq
\mathfrak{d}} \sharp \left (\mathfrak{a}/\mathfrak{d} \right
)^{-s}; ~ \mathrm{Re} (s) >0,
\]
with $\mathfrak{d}$ a fractional ideal of $\oo$. This series decomposes
as
\[
\zeta (\mathfrak{d},s)= \sum_{(\mathfrak{b})} \zeta
(\mathfrak{d},\mathfrak{b},s),
\]
where $(\mathfrak{b})$ means that the sum runs through a complete
system of representatives of the ideal class semigroup of $\oo$,
and the so-called \emph{partial} zeta functions:
\[
\zeta (\mathfrak{d}, \mathfrak{b}, s) := \sum_{\mathfrak{a}
\supseteq \mathfrak{d}}^{\mathfrak{a} \sim \mathfrak{b}} \sharp
\left ( \mathfrak{a} / \mathfrak{d}  \right )^{-s}, ~ ~
\mathrm{Re}(s) > 0.
\]
The sum runs through all the ideals $\mathfrak{a} \supseteq
\mathfrak{d}$ which are equivalent to $\mathfrak{b}$. Two ideals
$\mathfrak{a}$ and $\mathfrak{b}$ are said to be
\textbf{equivalent}, and we write $\mathfrak{a} \sim \mathfrak{b}$, if there exists $z \in K^{\ast}$ such that
$\mathfrak{a} = z^{-1} \mathfrak{b}$. 
\medskip

As $\sharp \left (\mathfrak{a} / \mathfrak{d} \right ) =
q^{\dim_{\ff} (\mathfrak{a}/\mathfrak{d})}$, we can rewrite the
partial zeta functions as power series in $t=q^{-s}$ (it is customary to use
the latin letter $Z$ when writing variable $t$ and  $\zeta$ in the case of using the variable $s$):
\[
Z(\mathfrak{d},\mathfrak{b},t):= \sum_{\mathfrak{a} \supseteq
\mathfrak{d}}^{\mathfrak{a} \sim \mathfrak{b}} t^{\dim_{\ff}
(\mathfrak{a} / \mathfrak{d})}
\]
and $Z(\mathfrak{d},t)=\sum_{(\mathfrak{b})}
Z(\mathfrak{d},\mathfrak{b},t)$, with $(\mathfrak{b})$ as usual.

\begin{defn}
For each fractional ideal $\mathfrak{a} \subseteq \mathfrak{\oo}$
we define the \textbf{degree} of $\mathfrak{a}$, denoted by $\deg
(\mathfrak{a})$, for the following two properties:
\begin{itemize}
  \item[i)] $\deg{\oo}=0$
  \item[ii)] $\dim_{\ff}{\mathfrak{a} /
  \mathfrak{b}}=\deg{\mathfrak{a}}-\deg{\mathfrak{b}}$ if $\mathfrak{a} \supseteq \mathfrak{b}$.
\end{itemize}
\end{defn}

%Dado $\mathfrak{d} \subseteq \mathfrak{a}$, al ser $\ago$ de la
%forma $\pi^{-\underline{n}} \mathfrak{b}$ con $\mathfrak{b}$ un
%ideal de $\oo$ que cumple $\mathfrak{b} \cdot \overline{\oo} =
%\overline{\oo}$:
%$$
%\pi^{-\underline{n}} \mathfrak{b} \supseteq \mathfrak{d}
%\Longleftrightarrow \pi^{\underline{n}} \in \mathfrak{b} :
%\mathfrak{d}.
%$$
%
%Definiendo, para cada ideal $I$:
%$$
%\Gamma (I) := \{ \underline{n} \in \mathbb{Z}^r \mid
%\pi^{\underline{n}} \in I \}
%$$
%se obtiene la partición
%$$
%Z(\mathfrak{d},t)= \sum_{\mathfrak{b} \cdot \overline{\oo} =
%\overline{\oo}} \left ( \sum_{\underline{n} \in \Gamma
%(\mathfrak{b} : \mathfrak{d} )} t^{\underline{n} \cdot
%\underline{d} + \deg (\mathfrak{b})-\deg(\mathfrak{d})} \right ).
%$$
%
%Esta partición no es canónica, pues $\Gamma (\mathfrak{a})$
%depende de la elección de los $\pi_i$. Por ello substituiremos los
%ideales $\pi^{-\underline{n}} \mathfrak{b}$, $\underline{n} \in
%\mathbb{Z}^r$, por $z^{-1} \mathfrak{b}$, $z \in K^{\ast}$, esto
%es, ideales \emph{equivalentes} a $\mathfrak{b}$. Así formamos
%funciones $\zeta$ locales parciales:
%$$
%\zeta (\mathfrak{d}, \mathfrak{b}, s) := \sum_{\mathfrak{a}
%\supseteq \mathfrak{d}}^{\mathfrak{a} \sim \mathfrak{b}} \sharp
%\left ( \mathfrak{a} / \mathfrak{d}  \right )^{-s}, ~ ~
%\mathrm{Re}(s) > 0
%$$

Next result is already proved in \cite{stohr}. We summarize here
his proof by the sack of completeness.

\begin{thm}[St\"ohr] \label{thm:31}
Let $U_{\mathfrak{b}}$ be the multiplicative group of elements $u
\in K$ such that $u \mathfrak{b}=\mathfrak{b}$. Then
\[
Z(\oo, \mathfrak{b},t)=\frac{1}{(U_{\mathfrak{b}} :U_{\oo})}
\sum_{\underline{n} \in S(\mathfrak{b})}
\varepsilon_{\underline{n}} (\mathfrak{b}) t^{\underline{n} \cdot
\underline{d} + \deg(\mathfrak{b})}
\]
with
\[
\varepsilon_{\underline{n}}(\mathfrak{b}):=
\frac{q^{\rho}}{q^{\rho}-1}
\sum_{\underline{i}=\underline{0}}^{(1,\ldots,1)} (-1)^{|i|}
q^{\underline{n} \cdot \underline{d} + \deg \left (\mathfrak{b}
\cap \mm^{\underline{n} + \underline{i}} \right )},
\]
where $\underline{i}=(i_1,\ldots,i_r) \in \{ 0,1\}^r $, $|i|=i_1+
\ldots + i_r$, and $\rho$ denotes the extension degree of the
residue field of $\oo$ over $\ff$.
\end{thm}

\dem

The ideals $\mathfrak{a}$ equivalent to $\mathfrak{b}$ such that
$\mathfrak{a} \supseteq \mathfrak{d}$, are those of the form
$z^{-1} \mathfrak{b}$, where $z$ varies on a complete system of
representatives of $(\mathfrak{b} : \mathfrak{d}) / \{ 0 \}$ under
the action of $U_{\mathfrak{b}}$. Then
\begin{eqnarray} 
\zeta (\mathfrak{d},\mathfrak{b},s) &=& \sum_{z \in (\mathfrak{b} :
\mathfrak{d}) \setminus \{ 0 \} / U_{\mathfrak{b}} } q^{-s
\dim_{\ff} (z^{-1}\mathfrak{b} / \mathfrak{d})} \nonumber \\
&=& (U_{\mathfrak{b}} : U_{\mathfrak{\oo}})^{-1} \sum_{z \in
(\mathfrak{b} : \mathfrak{d}) \setminus \{ 0 \} / U_{\oo} } q^{-s
\dim_{\ff} (z^{-1}\mathfrak{b} / \mathfrak{d})}. \nonumber
\end{eqnarray}

The index $(U_{\mathfrak{b}}:U_{\oo})$ is finite, since $U_{\oo}
\subseteq U_{\mathfrak{b}} \subseteq U_{\overline{\oo}}$ and
$(U_{\overline{\oo}} : U_{\oo})$ is also finite.
\medskip

We define, for every $z \in K$, the \emph{absolute value} of $z$ as
\[
|z|:=\left\{%
\begin{array}{ll}
    q^{\deg (z \oo)} = q^{ \deg (z \overline{\oo}) - \deg (\overline{\oo})} = q^{\sum_{i=1}^{r} d_i v_i (z)}~ \mathrm{if}~ z \ne 0 \\
    0 ~ ~ \mathrm{if}~ z=0.\\
\end{array}%
\right.
\]

\emph{Remark.} If $z \in K^{\ast}$, then $z \oo \subset z \ooo$,
and the multiplication by $z$ defines an automorphism of $\ooo$,
in such a way that the quotient vector spaces $z \ooo / z \oo$ and
$\ooo / \oo$ are isomorphic. Hence
\[
\deg \ooo = \dim_{\ff} \ooo /
\oo = \dim_{\ff} z \ooo / z \oo = \deg (z \ooo) - \deg (z \oo)
\]
and $\deg (z \oo) = \deg (z \overline{\oo}) - \deg
(\overline{\oo})$.
\medskip

Therefore

\begin{eqnarray}
\zeta (\mathfrak{d},\mathfrak{b},s) & = &
\frac{1}{(U_{\mathfrak{b}} : U_{\oo})} \sum_{z \in
(\mathfrak{b}:\mathfrak{d}) \setminus \{ 0 \} / U_{\oo}} q^{-s
\dim_{\ff}( z^{-1}\mathfrak{b}/ \mathfrak{d})} \nonumber \\
 & = & \frac{1}{(U_{\mathfrak{b}} : U_{\oo})}\sum_{z
 \in (\mathfrak{b}:\mathfrak{d}) \setminus \{ 0 \} / U_{\oo}} q^{-s \left ( \deg (z^{-1}\mathfrak{b}) - \deg (\mathfrak{d}) \right )} \nonumber \\
 & = & \frac{1}{(U_{\mathfrak{b}} : U_{\oo})} \sum_{z
 \in (\mathfrak{b}:\mathfrak{d}) \setminus \{ 0 \} / U_{\oo}} q^{-s \left ( \deg \mathfrak{b} - \deg (z \oo) - \deg \mathfrak{d} \right )} \nonumber \\
 & = & \frac{q^{-s(\deg(\mathfrak{b}) - \deg (\mathfrak{d})
)}}{(U_{\mathfrak{b}} : U_{\oo})} \sum_{z \in
(\mathfrak{b}:\mathfrak{d}) \setminus \{ 0 \} / U_{\oo}} q^{s \deg
(z \oo)} \nonumber \\
& = & \frac{q^{-s(\deg(\mathfrak{b}) - \deg (\mathfrak{d})
)}}{(U_{\mathfrak{b}} : U_{\oo})} \sum_{z \in (\mathfrak{b} :
\mathfrak{d}) / U_{\oo}} |z|^s . \nonumber
\end{eqnarray}

If $\mathfrak{d}= \oo$ we have
\begin{equation} \nonumber
\zeta(\oo,\mathfrak{b},s)= \frac{q^{-s\deg(\mathfrak{b})
}}{(U_{\mathfrak{b}} : U_{\oo})} \sum_{z \in \mathfrak{b} /
U_{\oo}} |z|^s.
\end{equation}

%Si usamos el hecho de que se puede asociar una medida de Haar
%$\widehat{\mu}$ al anillo localmente compacto de la compleción
%$\widehat{\oo}$ de $\oo$, la suma de (\ref{eq:27}) se interpreta
%como:
%$$
%\sum_{z \in (\mathfrak{b} : \mathfrak{d}) / U_{\oo}} |z|^s =
%\frac{q^{\rho}}{q^{\rho}-1}
%\int_{\widehat{(\mathfrak{b}:\mathfrak{d})}} |z|^{s-1} d
%\widehat{\mu} (z),
%$$
%donde el entero $\rho$ es el grado del cuerpo residual de $\oo$
%sobre el cuerpo base $\ff$.

Let $\mathcal{M}_0$ be the set of fractional ideals of $\oo$. It
can be extended to a system $\mathcal{M}$ of subsets $M \subset K$
which are obtained from $\mathcal{M}_0$ by addition of an element
of $K$. The system $\mathcal{M}$ is endowed with the operations
intersection and union.
\medskip

We associated a measure $\mu (M) \ge 0$ to each $M \in
\mathcal{M}$. This measure is uniquely determined by the following
three properties:
\begin{enumerate}
    \item $ \mu (\oo) =1$
    \item $ \mu (z+M) = \mu (M)$
    \item $ \mu(M \cup N) = \mu (M) + \mu (N) - \mu (M \cap N)$.
\end{enumerate}
\medskip

Let $\mathfrak{a}$ be an integral ideal of $\oo$. The local ring
$\oo$ is disjoint union of finitely many cosets of $\mathfrak{a}$
in the quotient ring $\oo / \mathfrak{a}$:
\[
\oo = (x_1 + \mathfrak{a}) \uplus \ldots \uplus (x_h +
\mathfrak{a}), ~ h = \sharp \left ( \oo / \mathfrak{a} \right ).
\]
Hence, $1 = \mu (\oo) = \mu (x_1 + \mathfrak{a}) + \ldots + \mu
(x_h + \mathfrak{a})= h \cdot \mu (\mathfrak{a})$, and
\[
\mu (\mathfrak{a})=\frac{1}{h} = \frac{1}{\sharp (\oo /
\mathfrak{a})} = q^{\deg (\mathfrak{a})}.
\]
Set $z \in K^{\ast}$. Consider the function
\[
M \mapsto \frac{\mu (zM)}{\mu (z \oo)}.
\]
Since it verifies the properties (1), (2) and (3) above and the
measure $\mu$ is unique, we obtain
\[
\mu (M)=\frac{\mu (zM)}{\mu (z \oo)} \Longleftrightarrow \mu (zM)
= \mu (M) \mu (z \oo).
\]

Now we want to compute $\mu (z \oo)$. If we write the element $z
\in K^{\ast}$ as a quotient $z = \frac{z_1}{z_2} \Leftrightarrow
z_2 \cdot z = z_1$, with $z_1,z_2 \in \oo$, then $\mu (z_2 \oo)
\mu (z \oo) = \mu(z \oo z_2) = \mu (z_1 \oo)$, and so
\[
\mu (z \oo)=\frac{\mu (z_1 \oo)}{\mu (z_2 \oo)}.
\]
Since $\mathfrak{a}$ is fractional, it holds that $z \mathfrak{a} \subset
\oo$, and we have
\[
\mu (z \mathfrak{a}) = \mu (z \oo) \mu (\mathfrak{a}).
\]
The measures $\mu (z \mathfrak{a})$ y $\mu (z \oo)$ are known,
therefore also $\mu (\mathfrak{a})=\frac{\mu (z \mathfrak{a})}{\mu
(z \oo)}$.
%$$
%\mu (\mathfrak{a}) = \frac{1}{\sharp \left ( \oo / \mathfrak{a}
%\right )}.
%$$
%Esto es,
%$$
%\mu (\mathfrak{a}) = q^{\deg (\mathfrak{a})}.
%$$
%Por la ecuación (\ref{eq:27}) sabemos que:
\medskip

For every $\underline{n} \in \mathbb{Z}^r$ we define
\[
\mathfrak{b}_{\underline{n}} := \mathfrak{b} \cap
\pi^{\underline{n}} U_{\overline{\oo}},
\]
where we note:
\begin{enumerate}
    \item $U_{\ooo}:= \{z \in \ooo \mid v_i(z)=0, ~ 1 \le i \le
    r\}$. If $\ooo \supset \mm_i$ for all $1 \le i \le r$, are the
    maximal ideals of $\ooo$, then $V_i:=\ooo_{\mm_i}$ is a discrete Manis valuation ring
    with asssociated Manis valuation $v_i$ for all $1 \le i \le r$. In this case, $\ooo = V_1 \cap
    \ldots \cap V_r$, and therefore $z \in U_{\ooo} \Leftrightarrow
    v_i(z)=0$ for $1 \le i \le r$.
    \item $\mathfrak{b}_{\underline{n}}= \{ z \in \mathfrak{b} \mid v_i (z)=n_i, ~ ~ 1 \le i \le
    r\}$, since $v_i(z)=n_i$ for all $1 \le i \le r$ if and only
    if
    $z=\pi^{\underline{n}}\varepsilon$, with $\varepsilon \in
    U_{\ooo}$.
    \item  $U_{\mathfrak{b}} \supset U_{\oo}$. %$U_{\mathfrak{b}}=\{
    %z \in K^{\times} \mid z \mathfrak{b}=\mathfrak{b}
    %\}$
\end{enumerate}

The elements $\mathfrak{b}_{\underline{n}}$ set a partition of
$\mathfrak{b} \setminus \{ 0\}$:
\[
\mathfrak{b} \setminus \{ 0 \} = \bigcup_{\underline{n} \ge
\underline{0}} \mathfrak{b}_{\underline{n}}.
\]
If $z \in \mathfrak{b}$, $\underline{n}=(v_1 (z), \ldots , v_r(z))
\in \mathbb{Z}^r$, then $\mathfrak{b}_{\underline{n}}=\{zu \in
\mathfrak{b} \mid u \in U_{\ooo} \}$.

Since $|z|=q^{-\underline{n} \cdot \underline{d}}$ when $z \in
\mathfrak{b}_{\underline{n}}$, putting $t=q^{-s}$ one has

\[
Z(\oo, \mathfrak{b},t) = \frac{t^{\deg
(\mathfrak{b})}}{(U_{\mathfrak{b}} : U_{\oo} )}
\sum_{\underline{n} \in S(\mathfrak{b})} \sharp \left (
\mathfrak{b}_{\underline{n}} / U_{\oo} \right )t^{\underline{n}
\cdot \underline{d}},
\]
where $S(\mathfrak{b}):= \{(v_1(z), \ldots , v_r(z)) \mid z \in \mathfrak{b} \setminus \{0\} \}$. The problem now is to compute the cardinality $\sharp \left (
\mathfrak{b}_{\underline{n}} / U_{\oo} \right )$. As each orbit of
the action of $U_{\oo}$ over $\mathfrak{b}_{\underline{n}}$ has
measure $\mu (U_{\oo} \pi^{\underline{n}}) = |\pi^{\underline{n}}|
\mu(U_{\oo}) = \mu (U_{\oo}) q^{-\underline{n} \cdot
\underline{d}}$, and since $\mu (U_{\oo})=1-q^{-\rho}$, we have
\[
\sharp \left ( \mathfrak{b}_{\underline{n}} / U_{\oo} \right ) =
\frac{\mu (\mathfrak{b}_{\underline{n}})}{\mu
(U_{\oo}\pi^{\underline{n}})}= \frac{q^{\underline{n} \cdot
\underline{d}}\mu (\mathfrak{b}_{\underline{n}}) }{1-q^{-\rho}}.
\]
Because of the fact that $\mathfrak{b}_{\underline{n}} = (\mathfrak{b} \cap
\mm^{\underline{n}}) \setminus \bigcup_{i=1}^{r} (\mathfrak{b}
\cap \mm^{\underline{n}}\mm_i)$, we obtain

\begin{eqnarray}
\mu (\mathfrak{b}_{\underline{n}}) & = & \mu (\mathfrak{b} \cap
\mm^{\underline{n}}) - \sum_{j=1}^{r} (-1)^{j-1} \sum_{i_1 <
\ldots < i_j} \mu (\mathfrak{b} \cap \mm^{\underline{n}} \mm_{i_1}
\cap \ldots \cap \mm^{\underline{n}} \mm_{i_j}) \nonumber
\\
& = & \sum_{j=0}^{r} (-1)^{j} \sum_{i_1 < \ldots < i_j} q^{\deg
(\mathfrak{b} \cap \mm^{\underline{n}} \mm_{i_1} \ldots
\mm_{i_j})} \nonumber
\\
& = & \sum_{|\underline{i}| = \underline{0}}^{(1, \ldots , 1)}
(-1)^{|\underline{i}|} q^{\deg (\mathfrak{b} \cap
\mm^{\underline{n}+\underline{i}})}, \nonumber
\end{eqnarray}
which proves the statement.
\qed
%Ello demuestra el resultado siguiente:

%\dem

%Ver \cite{stohr}, Theorem 3.1., p.180.
%\qed

%\begin{cor} \label{cor:115}
%Si $\underline{n}$ es suficientemente grande como para que
%$\mm^{\underline{n}} \subseteq \mathfrak{b}$, entonces
%$$
%\varepsilon_{\underline{n}} (\mathfrak{b}) =
%\frac{q^{\delta}}{1-q^{\rho}} \prod_{i=1}^{r} \left ( 1-q^{-d_i}
%\right ),
%$$
%donde $\delta = \dim_{\ff} (\overline{\oo} / \oo)$.
%\end{cor}
%
%\dem
%
%Si $\mm^{\underline{n}} \subseteq \mathfrak{b}$, entonces
%$$
%\varepsilon_{\underline{n}}(\mathfrak{b})=
%\frac{q^{\rho}}{q^{\rho}-1} \sum_{\underline{i}=
%\underline{0}}^{\underline{1}} (-1)^{|\underline{i}|}q^{\delta
%-\underline{i} \cdot \underline{d}} = \frac{q^{\rho +
%\delta}}{q^{\rho}-1} \prod_{i=1}^{r} \left ( 1-q^{-d_i} \right ).
%\qed
%$$
%
%Como aplicación de este corolario, poniendo
%$\mathfrak{b}=\overline{\oo}$ y $\underline{n}=\underline{0}$, se
%obtiene:
%$$
%(U_{\overline{\oo}}:U_{\oo})= \frac{q^{\delta}}{1-q^{-\rho}}
%\prod_{i=1}^{r} \left (1-q^{-d_i} \right ).
%$$

\begin{thm}[St\"ohr] 
For every ideal $\mathfrak{b}$ of $\oo$, the series $Z(\oo,
\mathfrak{b},t)$, and therefore also $Z(\oo,t)$, is a rational
function in $t$.
\end{thm}

\dem

See \cite{stohr}, (3.4), p. 182 and (3.10), p.186.\qed

\section{On the computation of the St\"ohr zeta function}

Along this section we shall consider totally rational points $P$,
that is, points having rational places over them (in the
normalization of the curve); in other words, we will assume $d_i=1$ for  every $i \in \{1, \ldots, r\}$.
\medskip

The relation between the St\"ohr zeta function and the value semigroup is given by a theorem of Z\'u\~niga (\cite{zu}, Th. 5.5., p.
86), which establishes that the partial zeta function
\[
Z(\oo,\oo,t)=\sum_{\nu=0}^{\infty} \sharp
\{\mathrm{integral~principal~ideals~ of~} \oo \mathrm{~of~
codimension~} \nu \} t^{\nu}
\]
is determined by the semigroup $S(\oo)=S$. From Theorem
\ref{thm:31} we obtain the same result, since
\[
Z(\oo,\oo,t)=\sum_{\underline{n} \in S}
\varepsilon_{\underline{n}} t^{n_1+\ldots+n_r},
\]
where
\[
\varepsilon_{\underline{n}}= \frac{q}{q-1} \sum_{j=0}^{r} (-1)^j
\sum_{i_1 < \ldots < i_j} q^{n_1+\ldots+n_r-\dim_{\ff} \left (\oo
/ \oo \cap \mm^{\underline{n}}\mm_{i_1} \cdot \ldots \cdot
\mm_{i_j} \right )}
\]
for each $\underline{n}=(n_1, \ldots , n_r) \in \mathbb{Z}_{+}^r$.
These coefficients can be written as depending on
\[
\ell \left (\underline{n} \right ):= \dim_{\ff} \left (\oo / \oo
\cap \mm^{\underline{n}} \right ).
\]
Such positive integers satisfy
\[
\ell (\underline{n}+\underline{e}_i) \le \ell (\underline{n}) +1
\]
for every $i=1, \ldots, r$ (with $\underline{e}_i$ the vector of value $1$ in the $i$th
coordinate and zero in the rest) and they determine the semigroup
of values in the sense of the following lemma.

\begin{lemma}
Preserving notations as above, for every $i=1, \ldots, r$, the following
statements are equivalent:
\begin{enumerate}
    \item $\ell (\underline{n}+\underline{e}_i) = \ell (\underline{n})
    +1$.
    \item $\oo \cap \mm^{\underline{n}}\mm_i$  $\subseteq  \ne$  $\oo \cap
    \mm^{\underline{n}}$.
    \item There exists $(\beta_1, \ldots,\beta_r) \in S$ with $\beta_i = n_i$ and $\beta_j \ge n_j$ for all $j \in \{1, \ldots , r\}$, $i \ne j$.
\end{enumerate}
\end{lemma}

We want to compute the dimensions $\ell (\underline{n})$. A first
result relates this dimension with the number of points of the
projected semigroups into $\Delta_J^k(\underline{n})$, with
$\sharp J = r-1$.

\begin{lemma}
Let $\underline{n}=(n_1, \ldots,n_r) \in \mathbb{Z}^r$ and  let
$m(\underline{n})$ be the number of maximal points of $S$ whose
$i$-th coordinate is less than $n_i$ for all $i \in \{1, \ldots, r
\}$. We have
\[
\ell (n_1, \ldots, n_r) \ge \ell (n_1,0, \ldots, 0) + \ell (0,
n_2, 0 , \ldots , 0) + \ldots + \ell (0, \ldots, 0, n_r) -
m(\underline{n}).
\]
\end{lemma}

\dem %Haremos la prueba por inducción en el número $r$ de ramas.

We know that $\ell (n_1, \ldots, n_r) = \ell (n_1, \ldots ,
n_{r-1},n_r-1)$ if and only if $(n_1, \ldots, n_r-1)$ is a gap.
This point is a gap since either $(0, \ldots, 0, n_r-1)$ is a gap
or the hyperplane of $\mathbb{Z}^r_{+}$ of equation $\{ X_r =
n_r-1 \}$ contains maximal points. Therefore, denoting by $\sharp
M(i)$ the number of maximal points in the hyperplane $\{X_r=n_r-i
\}$ for $1 \le i \le r$:
\begin{eqnarray}
  \ell (n_1, \ldots, n_r) & = & \ell (n_1, \ldots, n_{r-1},0)+n_r- \left
( \sharp \{ \mathrm{gaps~ of ~} S_r \} + \mathrm{~some~maximal~points} \right )  \nonumber \\
   & \ge & \ell (n_1, \ldots, n_{r-1})+n_r-
 \sharp \{ \mathrm{gaps~ of ~} S_r \} -\sharp M(1) \nonumber
\end{eqnarray}
\begin{eqnarray}
  \ell (n_1, \ldots, n_{r-1}) & = & \ell (n_1, \ldots, n_{r-2},0)+n_{r-1}- \left
( \sharp \{ \mathrm{gaps~ of ~} S_{r-1} \} + \mathrm{~some~max.~pt} \right )  \nonumber \\
   & \ge & \ell (n_1, \ldots, n_{r-2})+n_{r-1}- \sharp \{ \mathrm{gaps~ of ~} S_{r-1} \} -
 \sharp M(2) \nonumber
\end{eqnarray}
Taking into account that $n_{i}- \sharp \{ \mathrm{gaps~ of ~}
S_{i} \} = \ell (n_i)$ for all $1 \le i \le r$, recursively we
obtain the result.\qed

\medskip

Set $e_{i_1,\ldots,i_h}:= (\underbrace{1,1,
\ldots,1}_{h},\underbrace{0,0,\ldots,0}_{r-h})$ for $1 \le h \le
r$. We want to compute the values $d^h(\underline{n})$ of the
differences
\[
\ell (\underline{n} + e_{i_1,\ldots,i_h}) - \ell (\underline{n})
=: d^h(\underline{n})
\]
depending on which kind of point is $\underline{n}$. If
$\underline{n}$ belongs to the hyperplane of $\mathbb{Z}^r_{+}$
having coordinates $\{X_i \ge n_i + 1,~ \forall ~i=1,
\ldots, h-1;~ X_h = n_h \}$, then $\ell (\underline{n} + e_{i_1,
i_2, \ldots , i_h})- \ell (\underline{n} + e_{i_1, i_2, \ldots ,
i_{h-1}})=1$. In this sense we say that each difference $\ell
(\underline{n} + e_{i_1, \ldots , i_h})- \ell (\underline{n} +
e_{i_1, \ldots , i_{h-1}})$, for $1 \le h \le r$, corresponds to
some hyperplane of $\mathrm{Z}^r_{+}$. The set of such all $r$
hyperplanes, for all $1 \le h \le r$, covers $\mathbb{Z}^r_{+}$
and forms what is called the \textbf{\emph{principal} cover
configuration}. (There are different cover configurations, one for
each $e_{i_j}$, $1 \le j \le h $ being removed from the differences
\[
\ell
(\underline{n}+ e_{i_1 \ldots i_h}) - \ell (\underline{n} + e_{i_1
\ldots \widehat{i_j} \ldots i_h})
\]
for $ h \in \{1, \ldots , r\}$; the
principal one is that in which we remove $e_{i_h}$ in the
subtrahend). Moreover, each hyperplane composing the cover
configuration is called a \textbf{chart} of the configuration. For
every point $\underline{n} \in S$, it is said that a chart
\textbf{covers} $\underline{n}$ if $\underline{n}$ belongs to the
hyperplane of $\mathbb{Z}^r_{+}$ corresponding to this chart.
\medskip

%El conjunto de $r$ hiperplanos de tales características que vienen
%de diferencias del tipo $\ell (\underline{n} + e_{i_1, \ldots ,
%i_h})- \ell (\underline{n} + e_{i_1, \ldots , i_{h-1}})$, para $1
%\le h \le r$, recubren $\mathbb{Z}_{+}^{r}$ y forman una
%\textbf{configuración de recubrimiento}, que llamaremos
%\emph{principal}. Existen diferentes configuraciones de
%recubrimiento, según las diferencias $\ell (\underline{n}+ e_{i_1
%\ldots i_h}) - \ell (\underline{n} + e_{i_1 \ldots \widehat{i_j}
%\ldots i_h})$, para $1 \le j \le h$ y $1 \le h \le r$, y son todas
%ellas equivalentes.
%
%\medskip
%
%Cada paso $\ell (\underline{n} + e_{i_1 \ldots i_j}) - \ell
%(\underline{n} + e_{i_1 \ldots i_{j-1}})$, para $1 \le j \le r$,
%se corresponde como hemos visto con un cierto hiperplano de
%$\mathbb{Z}_{+}^r$, las $r$ diferencias recubren
%$\mathbb{Z}_{+}^r$ y por ello las denominaremos \textbf{cartas} de
%la configuración. Dado un punto $\underline{n} \in S$, diremos que
%una carta \emph{cubre} a $\underline{n}$ si $\underline{n}$
%pertenece al hiperplano de $\mathbb{Z}_{+}^r$ dado por la carta.
%\medskip
%
%Con toda esta terminología previa, estamos en condiciones de
%probar:

\begin{prop} \label{Prop:maximal}
If $\underline{n}=(n_1, \ldots , n_r)$ is a maximal point of kind
$x_2, x_3, \ldots, x_{r-1}$:
%$$
%\ell (\underline{n}+\underline{1}) = \ell (\underline{n}) +
%\sum_{i=2}^{r-1} x_i +1.
%$$
\[
d^h(\underline{n})= \left\{%
\begin{array}{ll}
    h, & \hbox{if}~ 1 \le h \le \sum_{j=2}^{r-1} x_j + 1 \\
    \sum_{j=2}^{r-1} x_j + 1, & \hbox{if}~ \sum_{j=2}^{r-1} x_j + 1 \le h \le r. \\
\end{array}%
\right.
\]
\end{prop}

\dem

Assume the principal cover configuration and its corresponding $r$
charts. If $\underline{n}$ is maximal of kind $x_2, \ldots,
x_{r-1}$, then $\sum_{j=2}^{r-1}x_j + 1$ are covered. So, the
differences $\ell(\underline{n} + e_{i_1, \ldots , i_h})-
\ell(\underline{n})=h$, because we have $h$ covered charts
whenever $h \le \sum x_j + 1$. For $\sum x_j + 1 \le h \le r$,
$d^h (\underline{n})$ is equal to the maximal number of charts
which are covered, i.e., $\sum x_j + 1$. \qed

\begin{remark} \label{obs:44}

Let $\underline{n} \in \mathbb{Z}_{+}^r$, $I=\{1, \ldots, r \}$,
$J \subset I$, $J \ne I$. Consider $\Delta_{I \setminus
J}(\underline{n})$. We call \emph{dimension} of $\Delta_{I
\setminus J}(\underline{n})$ to $\dim(\Delta_{I \setminus
J}(\underline{n}))= \sharp I - \sharp J$. For any $K \subseteq I$,
$K \ne I$, with $\sharp K =i$, any $J \subseteq I$, $J \ne I$ with
$\sharp J = j$ and $J^{\ast} \subseteq J$, $J^{\ast} \ne J$ with
$\sharp J^{\ast} = j-1$, the pivotage property $(S4)$ states
that, if $\Delta_{I \setminus K}(\underline{n}) \ne \emptyset$
with dimension $r-i$ and $\Delta_{I \setminus J}(\underline{n})
\ne \emptyset$ with dimension $r-j$, then $\Delta_{I \setminus
J^{\ast}}(\underline{n}) \ne \emptyset$ with dimension $r-(j-1)$,
and this for all $i \in \{1, \ldots, r-1 \}$ and all $j \in \{i+1,
\ldots, r \}$.

\end{remark}

\begin{prop}\label{Prop:nomaximal}
If $\underline{n}= (n_1 , \ldots, n_r)$ is a non-maximal point of
the semigroup, then we have
\[
d^h (\underline{n}) \le h, ~ ~ ~ 1 \le h \le r.
\]
\end{prop}

\dem

As $\underline{n}$ is non-maximal, there is $i \in I$ such that
$\Delta_i (\underline{n}) \ne \emptyset$. Since $\underline{n} \in
S$, by $(S4)$ there exists $I^{\ast}=I \setminus \{ i \}$ such
that $\Delta_{I^{\ast}}(\underline{n}) \ne \emptyset$, and by the
Remark (\ref{obs:44}), $\Delta_{I^{\ast}}(\underline{n}) \ne
\emptyset$ for any $I^{\ast}$ with $\sharp I^{\ast} = r-1$.
Applying recursively (\ref{obs:44}), then
$\Delta_{K}(\underline{n}) \ne \emptyset$ for all $K \subseteq I$,
$K \ne I$, and therefore all charts of any cover configuration are
covered, hence it is easily checked that the differences $\ell (\underline{n}+ e_{i_1,
\ldots , i_k}) - \ell (\underline{n}+e_{i_1, \ldots , i_{k-1}})$
are at most $1$ for all $k$, and thus $d^h (\underline{n})\le h$,
for $1 \le h \le r$.

\qed

Given a semigroup $S$ associated with $r$ branches, and using the
Theorem \ref{thm:31}, we want to give an upper bound for the
coefficients
\[
\varepsilon_{\underline{n}}= \frac{q}{q-1} \sum_{j=0}^{r} (-1)^j
\sum_{i_1< \ldots < i_j} q^{n_1 + \ldots + n_r - \dim_{\ff}\left(
\oo / \oo \cap \mathfrak{m}^{\underline{n}} \mathfrak{m}_{i_1}
\cdot \ldots \cdot \mathfrak{m}_{i_j} \right)}.
\]
of the Zeta Function
\[
Z(\oo,\oo,t)=\sum_{\underline{n} \in S}
\varepsilon_{\underline{n}} t^{n_1+ \ldots + n_r},
\]

In particular, as $P$ is rational, $\underline{n} \cdot
\underline{d} = n_1 + \ldots + n_r$, and moreover
\[
\ell (\underline{n}) = \dim_{\ff} \left ( \oo / \oo \cap
\mathfrak{m}^{\underline{n}} \right ) = \deg (\oo)- \deg (\oo \cap
\mathfrak{m}^{\underline{n}})= - \deg (\oo \cap
\mathfrak{m}^{\underline{n}}).
\]
Thus the coefficients are expressed as
\[
\varepsilon_{\underline{n}}= \frac{q}{q-1}
\sum_{\underline{i}=\underline{0}}^{\underline{1}} (-1)^{i_1+
\ldots + i_r} q^{\underline{n} \cdot \underline{d} + \deg \left (
\oo \cap \mathfrak{m}^{\underline{n}} \right )},
\]
where $\underline{i}:=(i_1, \ldots, i_r)$.

\begin{thm}

The zeta function of the semigroup associated with $r$ branches is
given by
\[
Z(\oo,\oo,t)=\sum_{\underline{n} \in S}
\varepsilon_{\underline{n}} t^{n_1+ \ldots + n_r},
\]
where, denoting by $s_i (n_i)$ the number of gaps of the semigroup
$S_i$ less than o equal to $n_i$, and by $m(\underline{n})$ the
number of maximal points with coordinates less than or equal to
the coordinates of $\underline{n}$, the coefficients admit the
following upper bounds:
\begin{enumerate}
    \item If $\underline{n}$ is maximal of kind $x_2, \ldots, x_{r-1}$:
    \[
\varepsilon_{\underline{n}} \le \frac{q^{s_1(n_1)+ \ldots + s_r
(n_r)+m(\underline{n})}}{1-q^{-1}} \left ( 1 +
q^{-\sum_{j=2}^{r-1}x_j} \left ( \sum_{h=1}^{r} (-1)^h
\binom{r}{h} q^{r-h-1} \right ) \right ).
    \]
    \item If $\underline{n}$ belongs to the semigroup, but it is
    non-maximal:
    \[
\varepsilon_{\underline{n}} \le \frac{q^{s_1(n_1)+ \ldots + s_r
(n_r)+m(\underline{n})}}{1-q^{-1}} \left ( 1 + \sum_{h=1}^{r}
(-1)^h \binom{r}{h} q^{h} \right ).
    \]
\end{enumerate}
\end{thm}

\dem

We have to compute upper bounds for the coefficients of the zeta
function $Z(\oo,\oo,t)$, from the Theorem \ref{thm:31}, that is
\[
\varepsilon_{\underline{n}} \overset{(\dag)}{=} \frac{q}{q-1} \sum_{j=0}^{r}
(-1)^j \sum_{i_1 < \ldots < i_j } q^{n_1 + \ldots + n_r - \ell
(\underline{n})}
\]
as $\underline{n} \cdot \underline{d} = n_1 + \ldots + n_r$
because $P$ is rational, and moreover $\ell (\underline{n}) =
\dim_{\ff} (\oo / \oo \cap \mathfrak{m}^{\underline{n}})= \deg
(\oo) - \deg (\oo \cap \mathfrak{m}^{\underline{n}})=-\deg(\oo
\cap \mathfrak{m}^{\underline{n}})$.

\medskip

If $\underline{n}$ is maximal of kind $x_2, \ldots , x_{r-1}$, then Proposition \ref{Prop:maximal} shows that
$\ell (n+e_{i_1})\overset{(\ast \ast)}{=}\ell (n) + \sum_{i=2}^{r-1} x_i + 1 - 1$. The first summand in $(\dag)$ (i.e., the one corresponding to $j=0$) is obviously $q^{-(\ell (\underline{n}))}$. The summands corresponding to $j=1$ are exactly
\[
- q ^{-\ell(n)} -q^{-\ell (n_1, \ldots , n_r + 1)} -q^{- \ell
(n_1, \ldots , n_{r-1}+1, n_r)} - \ldots - q^{- \ell (n_1 +1 ,
\ldots , n_r)} ,
\]
and because of $(\ast \ast)$ we have

\[
-q^{-\ell (n)- \sum x_i -1+1} - \ldots - q^{-\ell (n) - \sum x_i -
1 + 1 } = -q^{- \ell (n) - \sum x_i - 1 } \Big ( \sum_{i=1}^{{r \choose 1}} q
\Big )
\]
Now we compute the summands corresponding to $j=2$; in this case we have
$\ell (n + e_{i_1 i_2}) = \ell (n)+ \sum_{i=2}^{r-1}
x_i + 1 - 2$, and in the same manner as before we get

\[
q^{-\ell (n_1, \ldots ,n_{r-2}, n_{r-1}+1 n_r + 1)} +q^{- \ell
(n_1, \ldots ,n_{r-2} n_{r-1}, n_r+1)} + \ldots + q^{- \ell (n_1
+1 , n_2+1, n_3, \ldots , n_r)}
\]
having in this case ${r \choose 2}$ summands.
\medskip

Proceeding in this way for every $j$ we obtain

\begin{eqnarray}
\varepsilon_n & = & \frac{q}{q-1} q^{\sum_{j=1}^{r} n_j} \Big [
q^{-\ell (n)} - q^{- \ell (n)-\sum x_i -1} \left ( q^{r-1} {r
\choose 1 } \right )+ \nonumber \\ & + & q^{- \ell (n)-\sum x_i
-1} \left ( q^{r-2} {r \choose 2 } \right) + \ldots + (-1)^r q^{-
\ell (n)-\sum x_i -1} \left ( q^{r-r} {r \choose r }\right ) \Big
] \nonumber \\ & = & \frac{q^{\sum_{j=1}^{r} n_j}}{1-q^{-1}} \left
[ q^{- \ell (n)} + q^{-\ell (n)-\sum x_j -1} \left (
\sum_{h=1}^{r} (-1)^h {r \choose h} q^{r-h} \right ) \right ]
\nonumber
\\ & = & \frac{q^{n_1 + \ldots + n_r}}{1-q^{-1}} \left [ q^{- \ell (n)} + q^{-\ell (n)-\sum x_j}  \left ( \sum_{h=1}^{r} (-1)^h {r \choose h} q^{r-h-1} \right ) \right ] \nonumber \\ & = &
\frac{q^{n_1 + \ldots + n_r}}{1-q^{-1}} q^{- \ell (n)} \left [1+
q^{-\sum x_j } \left ( \sum_{h=1}^{r} (-1)^h {r \choose h}
q^{r-h-1} \right ) \right ] \nonumber
\\ & = & (\ast) \nonumber
%\varepsilon_n & = & \frac{q}{q-1} q^{\sum_{j=1}^{r} n_j} \left [
%q^{-\ell (n)} - q^{- \ell (n)-\sum x_i -1} \left ( q^{r-1} {r
%\choose 1 } \right )  \nonumber \\
%&  & q^{- \ell (n)-\sum x_i -1}
%\left ( q^{r-2} {r \choose 2 } \right) + \ldots + (-1)^r q^{- \ell
%(n)-\sum x_i -1} \left ( q^{r-r} {r \choose r }\right ) \right ]
%\nonumber \\
%  & = & \frac{q^{\sum_{j=1}^{r} n_j}}{1-q^{-1}} \left [ q^{- \ell (n)} + q^{-\ell (n)-\sum x_j -1} \left ( \sum_{h=1}^{r} (-1)^h {r \choose h} q^{r-h} \right ) \right ] \nonumber \\
%  & = & \frac{q^{n_1 + \ldots + n_r}}{1-q^{-1}} \left [ q^{- \ell (n)} + q^{-\ell (n)-\sum x_j}  \left ( \sum_{h=1}^{r} (-1)^h {r \choose h} q^{r-h-1} \right ) \right ] \nonumber \\
%  & = & \frac{q^{n_1 + \ldots + n_r}}{1-q^{-1}} q^{- \ell (n)} \left [1+ q^{-\sum x_j } \left ( \sum_{h=1}^{r} (-1)^h {r \choose h} q^{r-h-1} \right ) \right ] \nonumber
 % & = & d
\end{eqnarray}

Since $\ell (n) \ge \ell (n_1, 0, \ldots, 0) + \ldots + \ell (0,
\ldots, n_r) - m(n)$ and moreover $n_1-\ell (n_1)=n_1-\ell (n_1, 0, \ldots,
0) = \sharp \{ \mathrm{gaps~of~} S_1 ~\mathrm{less~than~} n_1 \}$,
if we denote $s_i (n_i):= n_i - \ell (n_i)$ for each $i$, one has that
\[
n_1 + \ldots + n_r - \ell (n_1, \ldots, n_r) \le s_1 (n_1) +
\ldots + s_r (n_r) + m(n),
\]
therefore
\begin{eqnarray}
(\ast) & \le & \frac{q^{s_1 (n_1) + \ldots + s_r (n_r) + m (n)}}{1
- q^{-1} } \left [ 1+ q^{-\sum x_j } \left ( \sum_{h=1}^{r} (-1)^h
{r \choose h} q^{r-h-1} \right ) \right ]. \nonumber
\end{eqnarray}

\medskip

If $\underline{n}$ is non-maximal, then we have $d^h (\underline{n})\le h$ by Proposition \ref{Prop:nomaximal}, hence
\begin{eqnarray}
  \varepsilon_{\underline{n}} & \le & \frac{q}{q-1} q^{\sum n_i} \left [ q^{-\ell (\underline{n})} - q^{-\ell (\underline{n})-1} {r \choose 1} + \ldots + q^{- \ell (\underline{n})-r} {r \choose r}   \right ] \nonumber \\
   & = & \frac{q}{q-1} q^{\sum n_i} q^{- \ell (\underline{n})} \left [ 1+ \sum_{h=1}^{r} (-1)^{h} q^{-h} {r \choose h}  \right ]
   \nonumber \\
   & \le & \frac{q^{s_1 (n_1) + \ldots + s_r (n_r) + m (\underline{n})}}{1-q^{-1}} \left [ 1+ \sum_{h=1}^{r} (-1)^{h} {r \choose h} q^{h}  \right]. \qed \nonumber
\end{eqnarray}

Notice that in the case of one and two branches we get the same formulae as those of St\"ohr (cf. \cite[Theorems 4.1, 4.3]{stohr}):

\begin{cor}
Let $r=1$, and write $L(\oo,\oo,t)=\sum_{i=0}^{f} n_i t^i$, for $f$ the conductor of $S$. Then we have $n_i = q^{s(i)}$ if $i \in S$ and $i-1 \notin S$, $n_i = - q^{s(i)}$ if $i \notin S$ and $i-1 \in S$, and $0$ otherwise.
\end{cor}

\begin{cor}
Let $r=2$. Then we have 
\begin{eqnarray}
\varepsilon_{(n_1,n_2)} &= & q^{s_1(n_1)+s_2(n_2)+m(n_1,n_2)}, \mathrm{~if~}   (n_1,n_2) \in S; \nonumber \\
\varepsilon_{(n_1,n_2)} &= & q^{s_1(n_1)+s_2(n_2)+m(n_1,n_2)} (\frac{q-1}{q}), \mathrm{~if~}  (n_1,n_2) \in S \setminus M, \nonumber
\end{eqnarray}
where $M$ is the set of maximal points of $S$.
\end{cor}

% ------------------------------------------------------------------------

%\subsection*{Acknowledgment}
%The results of this paper were obtained during my Ph.D. stayment
%at Paderborn University and are also contained in my
%thesis~\cite{serre} with the same title. I would like to express
%deep gratitude to my supervisor in Paderborn Prof. Dr. Karlheinz
%Kiyek and his wife Gisolde whose guidance and support were crucial
%for the successful completion of this project.

% ------------------------------------------------------------------------
%GATHER{Xbib.bib}   % For Gather Purpose Only
%GATHER{Paper.bbl}  % For Gather Purpose Only
%\bibliographystyle{amsplain}
%\bibliography{Xbib}

\end{document}